\documentclass[12pt]{amsart}
\usepackage{amssymb}
\usepackage{amscd}
\def\Ad{\mathop{\mathrm {Ad}}\nolimits}

\def\Aut{\mathop{\mathrm {Aut}}\nolimits}

\newtheorem{theorem}{Theorem}[section]
\newtheorem{proposition}[theorem]{Proposition}
\newtheorem{lemma}[theorem]{Lemma}
\newtheorem{remark}[theorem]{Remark}
\newtheorem{corollary}[theorem]{Corollary}
\newtheorem{definition}[theorem]{DEFINITION}
\begin{document}
\title{Quantum co-adjoint orbits of the real diamond group}
\author{ Nguyen Viet Hai}
\date {Version of December 25,1999}
\address{Haiphong Teacher's Training College ,Haiphong city,Vietnam}
\email{nguyen\_viet\_hai@yahoo.com}
\thanks{This work was supported in part by the Vietnam National  Foundation for Fundamental Science Research .}
\keywords{real diamond group,Moyal $\star$-product ,quantum half-plans,quantum hyperbolic cylinders,quantum hyperbolic paraboloids}
\maketitle

\begin{abstract}
We present explicit formulas for deformation quantization on the co-adjoint orbits of the real diamond Lie  group. From this we obtain quantum half-plans, quantum hyperbolic cylinders, quantum hyperbolic paraboloids via Fedosov deformation quantization and finally, the corresponding unitary representations of this group.
\end{abstract}

\section{Introduction}
 Let us first recall that it was Hermann Weyl(see \cite{weyl}), who introduced a mapping from classical observables (i.e. functions on the phase space  $\mathbb R^{2n}$) to quantum observables (i.e. normal operators in the Hilbert space $L^2(\mathbb R^n)$). The idea was to express functions on $\mathbb R^{2n}$ as Fourier transforms and then, by using the inverse Fourier transforms to correspond this correspondence to functions on characters, i.e. one-dimensional representations of the Heisenberg group, parameterized by the Planck constant $\hbar$ - and finally to present them as elements  in the corresponding infinite-dimensional representations of the Heisenberg group. This profound idea was later retrieved by Moyal, who have  seen that the symbols of the commutators or of the products of operators are of the the form of $sine$ (or $exponential$) functions of the bidifferential operators (of the Poisson brackets) of the corresponding symbols.

In the early 70's Berezin has treated the general mathematical definition of quantization as a kind of a functor from the category of classical mechanics to a certain category of associative algebras. About the same time as F. A. Berezin, M. Flato, M. G. Fronsdal, F. Bayen, A. Lichnerowicz and D. Sternheimer considered quantization as a deformation of the commutative products of classical observables into a noncommutative $\star $-products which are parameterized by the Planck constant $\hbar$ and satisfy the {\it correspondence principle}.
They systematically developed the notion of deformation quantization as a theory of $\star$-products and gave an independent formulation of quantum mechanics based on this notion (see\cite{rt}). 

It was proved by Gerstenhaber that a formal deformation quantization exists on an arbitrary symplectic manifold, see for example \cite{fedosov} for a detailed explaination. It is however formal and quite complecate in general. We would like to simplify it in some particular cases.
\par

From the orbit method, it is well-known that coadjoint orbits are homogeneous symplectic manifolds with respect to the natural Kirillov structure form on coadjoint orbits. A natural question is to associate in a reasonable way to these orbits some quantum objects, what could be called {\it quantum co-adjoint orbits}. In particular, in \cite{diephai1} and \cite{diephai2} we obtained  ``quantum half-planes'' and ``quantum punctured complex planes", associated with the affine transformation groups of the real or complex straightlines. In this paper we will therefore continue to realize the problem for  the real diamond Lie group. This group has a lot of nontrivial 2-dimensional coadjoint orbits, which are the  half-planes, the hyperbolic cylinders and the hyperbolic paraboloids. We should find out explicit formulas for each of these orbits. Our main result therefore is the fact that by using $\star$-product we can construct the corresponding quantum half-plans, quantum hyperbolic cylinders, quantum hyperbolic paraboloids and by an {\it exact computation} we can find out explicit $\star$-product formulas and then, the complete list of irreducible unitary representations of this group. It is useful to do here a remark that there is a general theory for exponetial and compact groups. But our consideration concerning with non-exponential and noncompact Lie group and associated $G$-homogeneous symplectic manifolds. 

\par Let us in  few words describe the structure od the paper.
We introduce some preliminary results in \S2. Then, the adapted chart and in particular, Hamiltonian functions in canonical coordinates of the co-adjoint orbit $\Omega_F$ are exposed in \S3. The operators $\hat{\ell}_A$ which define the representations of the real diamond Lie algebra  are constructed in \S4 and finally,  by exponentiating them, we obtain the corresponding unitary representations of the real diamond Lie group $\mathbb R \ltimes \mathbb H_3$.

\section{Preliminary results}

The so called real diamond Lie algebra is the 4-dimensional solvable Lie algebra   $\mathfrak g$ with basis X, Y, Z, T satisfying the following commutation relations:
$$[X,Y]=Z,[T,X]=-X,[T,Y]=Y,$$ $$[Z,X]=[Z,Y]=[T,Z]=0.$$  
These relations  show that this real diamond Lie algebra $\mathbb R\ltimes \mathfrak h_3$ is an extension of the one-dimensional Lie algebra $\mathbb R T$ by the Heisenberg algebra $\mathfrak h_3$ with basis  $X, Y, Z$, where the action of  T  on  Heisenberg algebra $\mathfrak h_3$  is defined by the matrix $$\ ad_T=\left(\begin{array}{ccc}-1&0&0\\0&1&0\\0&0&0 \end{array}\right).$$
We introduce the following notations. The real diamond Lie algebra is isomorphic to $\mathbb R^4$ as vector spaces. The coordinates in this standard basis is denote by $(a,b,c,d)$. We identify its dual vector space $\mathfrak g^*$ with $\mathbb R^4$ with the help of the dual basis $X^*,Y^*,Z^*,T^* $ and with the local coordinates as $(\alpha,\beta, \gamma,\delta)$. Thus, the general form of an element of $\mathfrak g $ is $U=aX+bY+cZ+dT$ and  the general form of an element of $\mathfrak g^* $ is   $F=\alpha{X^*} +\beta{Y^*}+\gamma{Z^*}+\delta{T^*}$. The co-adjoint action of $G=\mathbb R \ltimes \mathbb H_3$ on $\mathfrak g^*$ is given (see e.g. \cite{kirillov1}) by $$\langle K(g)F, U\rangle = \langle F, \Ad(g^{-1})U\rangle,\quad  \forall F \in \mathfrak g^*, g \in G \mbox{ and } U \in \mathfrak g=Lie(\mathbb R \ltimes \mathbb H_3).$$ Denote the co-adjoint orbit of $G$ in $\mathfrak g$, passing through $F$ by 
$$\Omega_F = K(G)F :=  \{K(g)F \vert g \in G \}.$$
By a direct computation one obtains (see \cite {dndiep}):
\begin {itemize}
\item Each point of the line $\alpha=\beta=\gamma =0$ is a 0-dimensional co-adjoint orbit
 \begin {equation} \Omega^1=\Omega_{(0,0,0,\delta)}.\end{equation}
\item The set $\alpha\ne 0,\beta=\gamma=0$ is union of 2-dimensional co-adjoint orbits ,which are just the\quad{\bf half-planes }
\begin{equation}\Omega^2=\{(x,0,0,t)\quad\vert\quad x, t\in \mathbb R , \alpha x>0\}.\end{equation} 
\item  The set $\alpha=\gamma=0,\beta\ne 0$ is a union of 2-dimensional co-adjoint orbits, which are \quad{\bf half-planes}
 \begin{equation}\Omega^3=\{(0,y,0,t)\quad\vert\quad y, t\in\mathbb R, \beta y >0 \}.\end{equation}
\item The set $\alpha\beta\ne 0,\gamma=0$ is decomposed into a family of 2-dimensional co-adjoint orbits, which are  {\bf hyperbolic\quad  cylinders }
\begin{equation}\Omega^4=\{(x,y,0,t)\quad\vert x,y,t\in\mathbb R\quad\&\quad \alpha x>0,\beta y>0,xy=\alpha\beta\}.\end{equation}
\item  The open set$ \gamma\ne 0$ is decomposed into a family of 2-dimensional co-adjoint orbits ,which are just the {\bf hyperbolic\quad  paraboloids}
 \begin{equation}\Omega^5=\{(x,y,\gamma,t)\quad\vert x,y,t \in\mathbb R \quad\&\quad xy- \alpha\beta=\gamma(t-\delta)\}.\end{equation}
\end{itemize}
Thus, the real diamond Lie algebra belongs to the class of $MD_4$-algebras , i.e. every K-orbit of the corresponding Lie group has dimension 0 or maximal (see \cite{dndiep}).

Let us consider now the problem of deformation quantization on  half-planes, hyperbolic\quad  cylinders, hyperbolic\quad  paraboloids. In order to do this, we shall construct   on  each of these orbits  a canonical Darboux coordinate system $(p,q)$ and a class of Hamiltonian functions in these coordinates.
 
\section{Hamiltonian functions in canonical coordinates of the  orbits $\Omega_F$}

Each element $A\in\mathfrak g$ can be considered as the restriction of the corresponding linear functional $\tilde A$ onto co-adjoint orbits, considered as a subset of $g^*$ ,$\tilde A(F)=\langle F,A\rangle$. It is well-known that this function is just the Hamiltonnian function, associated with the Hamiltonian vector field $\xi_A$, defined by the formula $$(\xi_Af)(x):=\frac {d}{dt}f(x\exp(tA))\vert_{t=0},\forall f\in C^\infty(\Omega _F).$$
It is well-known the relation $\xi_A(f)=\{\tilde A,f\},\forall f\in C^\infty(\Omega_F)$. Denote by $\psi$ the symplectomorphism from $\mathbb R^2$ onto $\Omega_F$
$$(p,q) \in \mathbb R^2 \mapsto \psi(p,q) \in \Omega_F,$$ we have:
\begin{proposition}
1. Hamiltonian function $ \tilde{A}$ in canonical coordinates $(p,q)$ of the orbit $\Omega_F$ is of the form $$\tilde{A}\circ\psi(p,q) =\cases
   $$dp+a\alpha e^{-q} ,$$ & \text{if \quad $ \Omega_F=\Omega^2$}\\
   $$dp+ b\beta e^q, $$ & \text{if \quad $ \Omega_F=\Omega^3$}\\
   $$ dp+a\alpha e^{-q}+b\beta e^q, $$&\text{if \quad $ \Omega_F=\Omega^4$}\\
  $$(d\pm b\gamma e^q)p\pm ae^{-q}\pm b(\alpha\beta -\gamma\delta)e^q+c\gamma, $$&\text{if \quad $ \Omega_F=\Omega^5$}\\
\endcases $$
2. In the canonical coordinates $(p,q)$ of the orbit $\Omega_F $, the Kirillov form $\omega$ is coincided with the standard form $ dp \wedge dq$.
\end{proposition}
{\it Proof}.
1. We adapt the diffeomorphism $\psi$ to  each of the following cases (for 2-dimensional co-adjoint orbits, only)
\begin{itemize}
\item With $\alpha\ne 0,\beta=\gamma=0$ 
$$(p,q) \in \mathbb R^2 \mapsto \psi(p,q)=(\alpha e^{-q},0,0,p) \in \Omega^2$$
Element $F\in \mathfrak g^*$ is of the form $F=\alpha{X^*} +\beta{Y^*}+\gamma{Z^*}+\delta{T^*}$, hence the value of the function $f_A=\tilde A$ on the element $A=aX+bY+cZ+dT$ is $\tilde A(F)=\langle F,A\rangle= $ $$\langle\alpha{X^*} +\beta{Y^*}+\gamma{Z^*}+\delta{T^*},aX+bY+cZ+dT \rangle=\alpha a+\beta b+\gamma c+\delta d .$$ It follows that \begin{equation}\tilde{A}\circ\psi(p,q) =a\alpha e^{-q} + dp,
\end{equation}
\item With $\alpha=\gamma=0,\beta\ne 0$,$$(p,q) \in \mathbb R^2 \mapsto \psi(p,q)=(0,\beta e^q,0,p) \in \Omega^3.$$$\tilde A(F)=\langle F,A\rangle =\alpha a+\beta b+\gamma c+\delta d . $   From this,
 \begin{equation}\tilde{A}\circ\psi(p,q) =b\beta e^q + dp\end{equation}
\item With $\alpha\beta\ne 0,\gamma=0$,
$$(p,q) \in \mathbb R^2 \mapsto \psi(p,q)=(\alpha e^{-q},\beta e^{q},0,p) \in \Omega^4.$$\begin{equation}\tilde{A}\circ\psi(p,q) =a\alpha e^{-q}+b\beta e^q + dp\end{equation}
\item At last, if  $ \gamma\ne 0$, we consider the orbit with  the first coordinate $ x>0$
$$(p,q) \in \mathbb R^2 \mapsto \psi(p,q)=(e^{-q},(\alpha\beta+\gamma p-\gamma\delta)e^q,\gamma,p) \in \Omega^5.$$ We have  \begin{equation}\tilde{A}\circ\psi(p,q) =ae^{-q}+b(\alpha\beta +\gamma p-\gamma\delta)e^q+c\gamma + dp=\end{equation}
$$ =(d+b\gamma e^q)p+ ae^{-q}+b(\alpha\beta-\gamma\delta)e^q+c\gamma.$$
The case $x<0$ is similarly treated:$$(p,q) \in \mathbb R^2 \mapsto \psi(p,q)=(-e^{-q},-(\alpha\beta+\gamma p-\gamma\delta)e^q,\gamma,p) \in \Omega^5.$$ \begin{equation}\tilde{A}\circ\psi(p,q) =-ae^{-q}-b(\alpha\beta +\gamma p-\gamma\delta)e^q+c\gamma + dp=\end{equation}
$$ =(d-b\gamma e^q)p- ae^{-q}-b(\alpha\beta-\gamma\delta)e^q+c\gamma.$$
\end{itemize}
2. We consider only  the following  case  (the rest are similar):\par
$$(p,q) \in \mathbb R^2 \mapsto \psi(p,q)=(e^{-q},(\alpha\beta+\gamma p-\gamma\delta)e^q,\gamma,p) \in \Omega^5.$$ $$\tilde{A}\circ\psi(p,q)  =(d+b\gamma e^q)p+ ae^{-q}+b(\alpha\beta-\gamma\delta)e^q+c\gamma.$$
In canonical Darboux coordinates $(p,q)$  ,$$F'=e^{-q}X^*+(\alpha\beta+\gamma p-\gamma\delta)e^q Y^*+\gamma Z^*+pT^*\quad\in\Omega^5,$$  and for $A=aX+bY+cZ+dT,\quad B=a'X+b'Y+c'Z+d'T$, we have    $\langle F',[A,B]\rangle=$
$$=\langle e^{-q}X^*+(\alpha\beta+\gamma p-\gamma\delta)e^q Y^*+\gamma Z^*+pT^*,(ad'-da')X+(db'-bd')Y+(ab'-ba')Z\rangle.$$It follows therefore that
\begin {equation}\langle F',[A,B]\rangle=(ad'-da')e^{-q}+(db'-bd')(\alpha\beta+\gamma p-\gamma\delta)e^q+\gamma(ab'-ba').
\end{equation}
On the other hand,   $$\xi_A(f)=\{\tilde{A},f\}=(d+b\gamma e^q)\frac{\partial f}{\partial q}-[-ae^{-q}+b(\alpha\beta+\gamma p-\gamma\delta)e^q]\frac{\partial f}{\partial p}$$
$$\xi_B(f)=\{\tilde{B},f\} =(d'+b'\gamma e^q)\frac{\partial f}{\partial q}-[-a'e^{-q}+b'(\alpha\beta+\gamma p-\gamma\delta)e^q]\frac{\partial f}{\partial p}.$$From this, consider two vector fields  $$\xi_A=(d+b\gamma e^q)\frac{\partial}{\partial q}-[-ae^{-q}+b(\alpha\beta+\gamma p-\gamma\delta)e^q]\frac{\partial}{\partial p},$$
$$\xi_B=(d'+b'\gamma e^q)\frac{\partial}{\partial q}-[-a'e^{-q}+b'(\alpha\beta+\gamma p-\gamma\delta)e^q]\frac{\partial}{\partial p}.$$

We have 
\begin{equation}\xi_A\otimes \xi_B =(d+b\gamma e^q)(d'+b'\gamma e^q) \frac{\partial}{\partial q}\otimes \frac {\partial}{\partial q}+\end{equation} $$+[(ad'-da')e^{-q}+(db'-d'b)(\alpha\beta+\gamma p-\gamma\delta)e^q]\frac{\partial}{\partial p}\otimes \frac {\partial}{\partial q}+$$ $$+[-ae^{-q}+b( \alpha\beta+\gamma p-\gamma\delta)e^q][-a'e^{-q}+b'( \alpha\beta+\gamma p-\gamma\delta)e^q]\frac{\partial}{\partial p}\otimes \frac {\partial}{\partial p}
$$
From (11) and (12) we conclude that in the canonical coordinates the Kirillov form is just the standard symplectic form $\omega = dp \wedge dq$.
The proposition is therefore proved.\hfill$\square$\par

\begin{definition}  Each chart  $\psi^{-1}$ on $\Omega_F$ which satisfy 1. and  2. of  proposition 3.1 is called an  adapted chart on $\Omega_F $.\end{definition}

In the next section we shall see that each adapted chart carries the Moyal $\star$-product from $\mathbb R^2$  onto  $\Omega_F$.

\section{Moyal $\star$-product and representations of $G=\mathbb R \ltimes \mathbb H_3.$}

Let us denote  by $\Lambda$ the  2-tensor associated with the Kirillov standard form $\omega =dp \wedge dq$  in canonical Darboux coordinates.  Let us consider the well-known Moyal $\star$-product of two smooth functions $u,v \in C^\infty(\mathbb R^{2n})$ (see e.g  \cite{arnalcortet1},\cite{diephai1}), defined by
$$u \star v = u.v + \sum_{r \geq 1} \frac{1}{r!}(\frac{1}{2i})^r P^r(u,v),$$ where
$$P^1(u,v)=\{u,v\}$$
$$P^r(u,v) := \Lambda^{i_1j_1}\Lambda^{i_2j_2}\dots \Lambda^{i_rj_r}\partial^r_{i_1i_2\dots i_r} u \partial^r_{j_1j_2\dots j_r}v,$$ with $$\partial^r_{i_1i_2\dots i_r} := \frac{\partial^r}{\partial x^{i_1}\dots \partial x^{i_r}};\quad x:= (p,q) = (p_1,\dots,p_n,q^1,\dots,q^n)$$ using multi-index notation. It is well-known that this series converges in the Schwartz distribution spaces $\mathcal S (\mathbb R^{2n})$. Furthermore, it was obtained the results (see e.g \cite{arnalcortet1}): If  $u,v \in \mathcal S (\mathbb R^{2n})$, then 
\begin{itemize}
\item $\overline u\star\overline v=\overline{v\star u}$
\item $\int(u\star v)(\xi)d\xi=\int uv d\xi $
\item $\ell_u : \quad  \mathcal S (\mathbb R^{2n})  \longrightarrow \mathcal S (\mathbb R^{2n}) $, defined by $\ell_u(v)= u\star v$ is continuous in $L^2(\mathbb R^{2n},d\xi)$ \quad and then can be extended to a bounded linear operator (still  denoted by $\ell_u$ ) on $L^2(\mathbb R^{2n},d\xi)$.
\end{itemize}
We apply this to the special case $n=1$, $x=(x^1,x^2)=(p,q)$  
\begin{proposition}
In the above mentioned canonical Darboux coordinates $(p,q)$ on the orbit $\Omega_F$, the Moyal $\star$-product satisfies the relation
$$i\tilde{A} \star i\tilde{B} - i\tilde{B} \star i\tilde{A} = i\widetilde{[A,B]}, \forall A, B \in \mathfrak g=Lie(\mathbb R \ltimes \mathbb H_3).$$
\end{proposition}
{\it Proof}.
We prove the proposition for the orbit $\Omega^5$, $\tilde{A} =(d+b\gamma e^q)p +ae^{-q}+b(\alpha\beta-\gamma\delta)e^q+c\gamma $ (the other cases are proved similar).
Consider the elements
 $A=aX+bY+cZ+dT,\quad B=a'X+b'Y+c'Z+d'T$, . Then as said above, the corresponding Hamiltonian functions are
$$\tilde{A} =(d+b\gamma e^q)p +ae^{-q}+b(\alpha\beta-\gamma\delta)e^q+c\gamma $$  $$\tilde{B} =(d'+b'\gamma e^q)p+a'e^{-q}+b'(\alpha\beta-\gamma\delta)e^q+c'\gamma $$
 It is easy then to see that
\begin {eqnarray*} 
&&P^0(\tilde{A},\tilde{B}) = \tilde{A}.\tilde{B}\\
&& P^1(\tilde{A},\tilde{B}) = \{ \tilde{A},\tilde{B}\} =  \partial_p\tilde{A}\partial_q \tilde{B} - \partial_q\tilde{A}\partial_p \tilde{B}=\\  
&& = (d+b\gamma e^q)[-a'e^{-q}+b'(\alpha\beta+\gamma p-\gamma\delta)e^q]-\\ &&-(d'+b'\gamma e^q)[-ae^{-q}+b(\alpha\beta+\gamma p-\gamma\delta)e^q]= \\
&&=[(ad'-da')e^{-q}+(db'-d'b)(\alpha\beta+\gamma p-\gamma\delta)e^q+(ab'-ba')\gamma]\\
&&P^2(\tilde{A},\tilde{B}) = \Lambda^{12}\Lambda^{12}\partial^2_{pp}\tilde{A} \partial^2_{qq}\tilde{B} +  \Lambda^{12}\Lambda^{21}\partial^2_{pq}\tilde{A} \partial^2_{qp}\tilde{B} + 
 \Lambda^{21}\Lambda^{12}\partial^2_{qp}\tilde{A} \partial^2_{pq}\tilde{B} +  \\&&+\Lambda^{21}\Lambda^{21}\partial^2_{qq}\tilde{A} \partial^2_{pp}\tilde{B} = -2bb'\gamma^2 e^{2q}\\
&& P^3(\tilde{A},\tilde{B}) = \Lambda^{12}\Lambda^{12}\Lambda^{12}\partial^3_{ppp}\tilde{A} \partial^3_{qqq}\tilde{B} +  \Lambda^{12}\Lambda^{12}\Lambda^{21}\partial^3_{ppq}\tilde{A} \partial^3_{qqp}\tilde{B} + \\&& +
\Lambda^{12} \Lambda^{21}\Lambda^{12}\partial^3_{pqp}\tilde{A} \partial^3_{qpq}\tilde{B} +  \Lambda^{21}\Lambda^{12}\Lambda^{12}\partial^3_{qpp}\tilde{A} \partial^3_{pqq}\tilde{B}+ \\&& + \Lambda^{21}\Lambda^{21}\Lambda^{12}\partial^3_{qqp}\tilde{A} \partial^3_{ppq}\tilde{B}+ \Lambda^{21}\Lambda^{12}\Lambda^{21}\partial^3_{qpq}\tilde{A} \partial^3_{pqp}\tilde{B}+\\&& + \Lambda^{12}\Lambda^{21}\Lambda^{21}\partial^3_{pqq}\tilde{A} \partial^3_{qpp}\tilde{B}+\Lambda^{21}\Lambda^{21}\Lambda^{21}\partial^3_{qqq}\tilde{A} \partial^3_{ppp}\tilde{B}=0 .  
\end{eqnarray*}
By analogy, we have $$P^k(\tilde{A},\tilde{B}) = 0, \forall k \geq 4.$$ Thus,
$$ i\tilde{A} \star i\tilde{B} - i\tilde{B} \star i\tilde{A} = \frac{1}{2i}[P^1(i\tilde{A} , i\tilde{B}) - P^1(i\tilde{B}, i\tilde{A})]$$ $$ = i[(ad'-da')e^{-q}+(db'-d'b)(\alpha\beta+\gamma p-\gamma\delta)e^q+(ab'-a'b)\gamma].$$
On the other hand, as  $$[A,B] = [aX+bY+cZ+dT,a'X+b'Y+c'Z+d'T]$$ $$=(ad'-da')X+(db'-d'b)Y+(ab'-a'b)Z $$ we obtain  $$ i[(ad'-da')e^{-q}+(db'-d'b)(\alpha\beta+\gamma p-\gamma\delta)e^q+(ab'-a'b)\gamma] $$ $$=i\widetilde{[A,B]} =  i\tilde{A} \star i\tilde{B} - i\tilde{B} \star i\tilde{A} .$$ The proposition is hence proved.\hfill$\square$\par

Consequently, to each adapted chart, we associate a $G$-covariant $\star$-product.Then  there exists a representation $\tau$ of $G$ in $\Aut N[[\nu]]$ ,(see \cite{gutt})  such that (here  $\nu=\frac i 2$):
$$\tau(g)(u \star v) = \tau(g)u \star \tau(g)v.$$
For each $A \in  Lie(\mathbb R\ltimes\mathbb H_3)$, the corresponding Hamiltonian function is $\tilde{A}$ and we can put $\ell_A(u) = i\tilde{A} \star u$,$u\in L^2(\mathbb R^2, \frac{dpdq}{2\pi})^\infty$. It is then continuated to the whole space $L^2(\mathbb R^2, \frac{dpdq}{2\pi})$. Because of the relation in Proposition (4.1), we have
\begin{corollary}
\begin{equation}
\ell_{[A,B]}=\ell_A\star\ell_B -\ell_B\star\ell_A:=[\ell_A,\ell_B]^\star
\end{equation}
\end {corollary}

This implies that  the correspondence $A \in  Lie(\mathbb R\ltimes\mathbb H_3) \mapsto \ell_A = i\tilde{A} \star .$ is a representation of the Lie algebra $ Lie(\mathbb R\ltimes\mathbb H_3)$ on the space $N[[\frac{i}{2}]]$ of formal power series in the parameter $\nu = \frac{i}{2}$(i.e $\hbar =1$) with coefficients in $N = C^\infty(M,\mathbb R)$ \cite{gutt}. 

\par
Let us denote by $\mathcal F_p(f)$ the partial Fourier transform  of the function $f$ from the variable $p$ to the variable $x$(see e.g\cite{meisevogt}), i.e.
$$\mathcal F_p(f)(x,q) := \frac{1}{\sqrt{2\pi}}\int_{\mathbb R} e^{-ipx} f(p,q)dp.$$ Let us denote by $ \mathcal F^{-1}_p(f) (p,q)$ the inverse Fourier transform. 

\begin{lemma}

1.$\partial_p \mathcal F^{-1}_p(f) = i \mathcal F^{-1}_p(x.f)$\par
2.$\mathcal F_p(p.v) = i \partial_x\mathcal F_p(v)  $\par
3.$\forall k\geq 2$ ,then
$P^k(\tilde{A},\mathcal F^{-1}_p(f))=$
$$=\cases $$a\alpha e^{-q}\partial^k_{p\dots p}\mathcal F^{-1}_p(f) $$ &\text{if $\tilde {A}$ is defined by (6) }\\
$$(-1)^k b\beta e^q\partial^k_{p\dots p}\mathcal F^{-1}_p(f)$$ &\text{if $\tilde {A}$ is defined by (7)}\\
$$[a\alpha e^{-q}+(-1)^k b\beta e^q] \partial^k_{p\dots p}\mathcal F^{-1}_p(f)$$ &\text{if $\tilde {A}$ is defined by  (8)} \\
$$(-1)^{k-1}k.b\gamma e^q \partial^k_{qp\dots p}\mathcal F_p^{-1}(f)+$$\\
$$+ [a e^{-q}+(-1)^k b(\alpha\beta+\gamma p-\gamma\delta)e^q]\partial^k_{p\dots p}\mathcal F_p^{-1}(f) $$&\text{if $\tilde {A}$ is defined by  (9) }\\
\endcases$$
\end{lemma}
{\it Proof}. The first two formulas are well-known from theory of Fourier transforms.

\par Let us prove 3. Remark that $\Lambda = \left(\begin{array}{cc} 0 & -1\\ 1 & 0 \end{array}\right)$ in the standard symplectic Darboux coordinates $(p,q)$ on the orbit $\Omega_F$, then 
\begin {itemize}
\item If $ \tilde A=a\alpha e^{-q}+dp $
\begin {eqnarray*} 
&& P^2(\tilde{A},\mathcal F^{-1}_p(f)) = \Lambda^{12}\Lambda^{12}\partial^2_{pp}\tilde{A} \partial^2_{qq}\mathcal F^{-1}_p(f))  +  \Lambda^{12}\Lambda^{21}\partial^2_{pq}\tilde{A} \partial^2_{qp}\mathcal F^{-1}_p(f))  + \\
&&\Lambda^{21}\Lambda^{12}\partial^2_{qp}\tilde{A} \partial^2_{pq}\mathcal F^{-1}_p(f))   +\Lambda^{21}\Lambda^{21}\partial^2_{qq}\tilde{A} \partial^2_{pp}\mathcal F^{-1}_p(f))=a\alpha e^{-q}\partial^2_{pp}\mathcal F^{-1}_p(f)=\\
&&P^3(\tilde A,\mathcal F^{-1}_p(f))=(-1)^6a\alpha e^{-q}\partial^3_{ppp}\mathcal F^{-1}_p(f)=a\alpha e^{-q}\partial^3_{ppp}\mathcal F^{-1}_p(f)
\end{eqnarray*}
and $P^k(\tilde A,\mathcal F^{-1}_p(f))=a\alpha e^{-q}\partial^k_{p\dots p}\mathcal F^{-1}_p(f)\quad\forall k \geq 4 ,$
\item If  $\tilde A=b\beta e^q+dp.$
$$P^k(\tilde A,\mathcal F^{-1}_p(f))=(-1)^k b\beta e^q\partial^k_{p\dots p}\mathcal F^{-1}_p(f)$$with  $\quad\forall k \geq 2 $
\item If $ \tilde{A}=a\alpha e^{-q}+b\beta e^q+dp,$
\begin {eqnarray*} 
&&P^2(\tilde{A},\mathcal F^{-1}_p(f)) = \Lambda^{12}\Lambda^{12}\partial^2_{pp}\tilde{A} \partial^2_{qq}\mathcal F^{-1}_p(f))  +  \Lambda^{12}\Lambda^{21}\partial^2_{pq}\tilde{A} \partial^2_{qp}\mathcal F^{-1}_p(f))  + \\
&&\Lambda^{21}\Lambda^{12}\partial^2_{qp}\tilde{A} \partial^2_{pq}\mathcal F^{-1}_p(f))   +\Lambda^{21}\Lambda^{21}\partial^2_{qq}\tilde{A} \partial^2_{pp}\mathcal F^{-1}_p(f))=\\
&&  =[a\alpha e^{-q}+(-1)^2b\beta e^q]\partial^2_{pp}\mathcal F^{-1}_p(f)\\
&&P^3(\tilde{A},\mathcal F^{-1}_p(f)) = [a\alpha e^{-q}+(-1)^3b\beta e^q]\partial^3_{ppp}\mathcal F^{-1}_p(f)
\end{eqnarray*}
By analogy we have $$P^k(\tilde{A}, \mathcal F^{-1}_p(f) )
 = [a\alpha e^{-q}+(-1)^k b\beta e^q] \partial^k_{p\dots p}\mathcal F^{-1}_p(f)),\quad \forall\quad k \geq 3.$$
\item If\quad $ \tilde{A}=(d+b\gamma e^q)p +ae^{-q}+b(\alpha\beta -\gamma\delta)e^q +c\gamma, $
\begin {eqnarray*} 
&&P^2(\tilde{A},\mathcal F^{-1}_p(f)) = \Lambda^{12}\Lambda^{12}\partial^2_{pp}\tilde{A} \partial^2_{qq}\mathcal F^{-1}_p(f))  +  \Lambda^{12}\Lambda^{21}\partial^2_{pq}\tilde{A} \partial^2_{qp}\mathcal F^{-1}_p(f))  + \\
&&\Lambda^{21}\Lambda^{12}\partial^2_{qp}\tilde{A} \partial^2_{pq}\mathcal F^{-1}_p(f))   +\Lambda^{21}\Lambda^{21}\partial^2_{qq}\tilde{A} \partial^2_{pp}\mathcal F^{-1}_p(f)) =\\
&&=(-1) 2. b\gamma e^q \partial_{qp}\mathcal F_p^{-1}(f)+[ae^{-q}+(-1)^2b(\alpha\beta+\gamma p-\gamma\delta)e^q]\partial^2_{pp}\mathcal F_p^{-1}(f).\\
&&P^3(\tilde{A},\mathcal F^{-1}_p(f))=(-1)^2.3b\gamma e^q\partial_{qpp}\mathcal F_p^{-1}(f)+\\&&+[ae^{-q}+(-1)^3b(\alpha\beta+\gamma p-\gamma\delta)e^q]\partial^3_{ppp}\mathcal F_p^{-1}(f).
\end{eqnarray*}

From this we also obtain :
\begin {eqnarray*} 
&&P^k(\tilde{A}, \mathcal F^{-1}_p(f) )=\\
&&(-1)^{k-1}.k.b\gamma e^q\partial^k_{qp\dots p}\mathcal F_p^{-1}(f)
 +\\&&+[ae^{-q}+(-1)^k b(\alpha\beta+\gamma p-\gamma\delta)e^q ]\partial^k_{p\dots p}\mathcal F^{-1}_p(f).\quad\forall\quad k\geq 3 
\end{eqnarray*}
\end{itemize}
The lemma is therefore proved.\hfill$\square$

\par
We study now the convergence of the formal power series. In order to do this, we look at the $\star$-product of $i\tilde{A}$ as the $\star$-product of symbols and define the differential operators corresponding to $i\tilde{A}$.  

\begin{theorem}
For each $A \in Lie(\mathbb R \ltimes \mathbb H_3)$ and for each compactly supported $C^\infty$ function $f \in C^\infty_0(\mathbb R^2)$, putting $\hat{\ell}_A(f) := \mathcal F_p \circ \ell_A \circ \mathcal F^{-1}_p(f)$,we have 
$$\hat{\ell}_A(f)  = \cases $$[d(\frac{1}{2}\partial_q-\partial_x)+ia\alpha e^{-(q-\frac{x}{2})}]f $$ &\text{if $\tilde {A}$ is defined by (6) }\\
$$[d(\frac{1}{2}\partial_q-\partial_x)+ib\beta e^{(q-\frac{x}{2})}]f $$ &\text{if $\tilde {A}$ is defined by  (7)}\\
$$[d(\frac{1}{2}\partial_q-\partial_x)+i(a\alpha e^{-(q-\frac{x}{2})}+ b\beta e^{(q-\frac{x}{2})})]f $$ &\text{if $\tilde {A}$ is defined by  (8)}$$ \\
$$[(d+ b\gamma e^{q-\frac{x}{2}})(\frac{1}{2}\partial_q-\partial_x)]f+$$\\ $$+i[ae^{-(q-\frac{x}{2})}+b(\alpha\beta-\gamma\delta)e^{q-\frac{x}{2}}+c\gamma]f $$&\text{if $\tilde {A}$ is defined by  (9)$$ }\\
$$[(d- b\gamma e^{q-\frac{x}{2}})(\frac{1}{2}\partial_q-\partial_x)]f+$$\\ $$+i[-ae^{-(q-\frac{x}{2})}-b(\alpha\beta-\gamma\delta)e^{q-\frac{x}{2}}+c\gamma]f $$&\text{if $\tilde {A}$ is defined by  (10)$$ }\\
\endcases$$
\end{theorem}
{\it Proof}. Applying Lemma (4.3),we have :
\begin{enumerate}
\item  If $\tilde A=a\alpha e^{-q}+dp$ then $$\hat{\ell}_A(f) := \mathcal F_p \circ \ell_A \circ \mathcal F^{-1}_p(f)  = \mathcal F_p(i\tilde{A}\star \mathcal F^{-1}_p(f)) = i\mathcal F_p\left(\sum_{r \geq 0} \left(\frac{1}{2i}\right)^r\frac{1}{r!} P^r(\tilde{A},\mathcal F^{-1}_p(f))\right)=$$
$$ \begin{array}{cl}&=i\mathcal F_p\Big\{(a\alpha e^{-q}+dp)\mathcal F_p^{-1}(f)+\frac{1}{1!}\frac{1}{2i}[d\partial_q\mathcal F_p^{-1}(f)+a\alpha e^{-q}\partial_p\mathcal F_p^{-1}(f)]+\\&+\frac{1}{2!}(\frac{1}{2i})^2.a\alpha e^{-q}\partial_pp^2\mathcal F_p^{-1}(f)+
\dots +\frac{1}{r!}(\frac{1}{2i})^r  a\alpha e^{-q}\partial^r_{p\dots p}\mathcal F_p^{-1}(f)+\dots \Big\}=\\& =i\Big\{a\alpha e^{-q}f+d\mathcal F_p(p.\mathcal F_p^{-1}(f))+\frac{1}{1!}\frac{1}{2i}[d\partial_qf +a\alpha e^{-q}\mathcal F_p(\partial_p\mathcal F_p^{-1}(f))]+\\& +\frac{1}{2!}(\frac{1}{2i})^2.a\alpha e^{-q}\mathcal F_p(\partial^2_{pp}\mathcal F_p^{-1}(f))+
+\frac{1}{3!}(\frac{1}{2i})^3.a\alpha e^{-q}\mathcal F_p(\partial^3_{ppp}\mathcal F_p^{-1}(f))+\dots\\&+\frac{1}{r!}(\frac{1}{2i})^r.a\alpha e^{-q}\mathcal F_p(\partial^r_{p\dots p}\mathcal F_p^{-1}(f))+\dots \Big\}=\\&
=d(\frac 12\partial_q-\partial_x)f+i a\alpha e^{-q}[1+\frac{x}{2}+\frac{1}{2!}(\frac{x}{2})^2+\dots +\frac{1}{r!}(\frac{1}{x})^r+\dots ]f =\\&= d(\frac 12\partial_q-\partial_x)f+i a\alpha e^{-q}e^{\frac{x}{2}}f=d(\frac 12\partial_q-\partial_x)f+ia\alpha e^{-(q-\frac{x}{2})}f\end{array}$$
\item  If $\tilde A=b\beta e^q +dp$ then    $$\hat\ell_A(f)=d(\frac 12\partial_q-\partial_x)f+i b\beta e^{q-\frac{x}{2}}f$$
\item  For each $\tilde A=a\alpha e^{-q}+b\beta e^q +dp ,$   we have:
$$ \begin{array}{cl}&\hat\ell_A=i\mathcal F_p\Big\{(a\alpha e^{-q}+b\beta e^q+dp)\mathcal F_p^{-1}(f)+\frac{1}{2i}[d\partial_q\mathcal F_p^{-1}(f)-(-a\alpha e^{-q}+\\&+ b\beta e^q)\partial_p\mathcal F_p^{-1}(f)]+\frac{1}{2!}(\frac{1}{2i})^2[a\alpha e^{-q}+(-1)^2b\beta e^q]\partial^2_{pp}\mathcal F_p^{-1}(f)+\dots +\\&+\frac{1}{r!}(\frac{1}{2i})^r[a\alpha e^{-q}+(-1)^r b\beta e^q]
 \partial^r_{p\dots p}\mathcal F_p^{-1}(f)+\dots\Big\}\\
&=ia\alpha e^{-q}.f+id\mathcal F_p\big(p.\mathcal F_p^{-1}(f)\big)+ib\beta e^q f+\frac{1}{2}d\partial_q f+\frac{1}{2}a\alpha e^{-q}\mathcal F_p\big(\partial_p\mathcal F_p^{-1}(f)\big)-\\
&-\frac{1}{2}b\beta e^q\mathcal F_p\big(\partial_p\mathcal F_p^{-1}(f)\big)+\dots i\frac{1}{r!}(\frac{1}{2i})^r a\alpha e^{-q}\mathcal F_p\big(\partial^r_{p\dots p}\mathcal F_p^{-1}(f)\big)+\\
&+i\frac{1}{r!}(\frac{-1}{2i})^r b\beta e^q\mathcal F_p\big(\partial^r_{p\dots p}\mathcal F_p^{-1}(f)\big)+\dots
=d(\frac{1}{2}\partial_q-\partial_x)+ia\alpha e^{-q}[1+\frac{x}{2}+\dots\\& +\frac{1}{r!}(\frac{x}{2})^r+\dots]+ib\beta e^q[1+(\frac{-x}{2})+\dots +\frac{1}{r!}(\frac{-x}{2})^r+\dots ]\\&
=d(\frac{1}{2}\partial_q-\partial_x)+i[a\alpha e^{-(q-\frac{x}{2})}+b\beta e^{q-\frac{x}{2}}]f.\end{array}$$
\item  For each $\tilde {A}$ is as in  (9) ,
 remark that $$ P^0(\tilde{A},\mathcal F^{-1}_p(f))=\tilde{A}.\mathcal F^{-1}_p(f);$$ $$P^1(\tilde{A},\mathcal F^{-1}_p(f)) = \{\tilde{A},\mathcal F^{-1}_p(f) \} =$$
$$(d+b\gamma e^q)\partial_q\mathcal F_p^{-1}(f)-[-ae^{-q}+b(\alpha\beta+\gamma p -\gamma\delta)e^q]\partial_p\mathcal F_p^{-1}(f)$$ and applying Lemma (4.3), we obtain:\quad 
$$\begin{array}{cl}& \hat{\ell}_A(f) = i\Big\{\mathcal F_p\Big([dp+ae^{-q}+b(\alpha\beta+\gamma p-\gamma\delta)e^{q}+c\gamma]\mathcal F_p^{-1}(f)\Big)+\\&+ \frac{1}{2i}\frac{1}{1!}\mathcal F_p\Big([d+b\gamma e^q]\partial_q\mathcal F^{-1}_p(f)-[-ae^{-q}+b(\alpha\beta+\gamma p -\gamma\delta)e^q]\partial_p\mathcal F_p^{-1}(f)\Big)+\\& +(\frac{1}{2i})^2\frac{1}{2!}\mathcal F_p\Big(-2b\gamma e^ {q}\partial^2_{pq}\mathcal F^{-1}_p(f)+ [ae^{-q}+b(\alpha\beta+\gamma p-\gamma\delta)e^q]\partial^2_{pp}(\mathcal F^{-1}_p(f)\Big)+\dots\\& +(\frac{1}{2i})^r\frac{1}{r!}\mathcal F_p\Big((-1)^{r-1}rb\gamma e^q\partial^r _{p \dots pq}\mathcal F^{-1}_p(f)+(-1)^r[(-1)^rae^{-q}+b(\alpha\beta+\gamma p-\gamma\delta)e^q]\times \\& \times \partial^r_{p\dots p}\mathcal F^{-1}_p(f)\Big)+\dots \Big\}=\end{array}$$
$$\begin{array}{cl}&=i\Big\{ae^{-q}f+b(\alpha\beta-\gamma\delta)e^qf+d\mathcal F_p\big(p\mathcal F^{-1}_p(f)\big) +b\gamma e^q\mathcal F_p\big(p\mathcal F^{-1}_p(f)\big)\\& +\frac{1}{2i}\frac{1}{1!} (d+b\gamma e^q)\partial_qf -\frac{1}{2i}[-ae^{-q}ixf  +b(\alpha\beta-\gamma\delta)e^qixf+b\gamma e^q\mathcal F_p\big(p\mathcal F^{-1}_p(xf)\big)] +\\&+(\frac{1}{2i})^2\frac{1}{2!}(-2b\gamma e^q)\mathcal F_p\left(\partial_{pq}^2\mathcal F^{-1}_p(f)\right)+(\frac {1}{2i})^2\frac{1}{2!}[ae^{-q}(ix)^2f+ b(\alpha\beta-\gamma\delta)e^q(ix)^2f+\\&+b\gamma e^q \mathcal F_p\left(pi^2\mathcal F^{-1}_p(x^2f)\right)]+\dots +(\frac {1}{2i})^r\frac{1}{r!}(-1)^{r-1}rb\gamma e^q\partial^r_{p\dots pq}\mathcal F^{-1}_p(f) \\&+(\frac{1}{2i})^r\frac{1}{r!}[ae^{-q}(ix)^rf+(-1)^rb(\alpha\beta-\gamma\delta)e^q(ix)^rf+b\gamma e^q\mathcal F_p\left(p(ix)^r\mathcal F^{-1}_p(f)\right)]+\dots \Big\}\\&=i[ae^{-q}(1+\frac{1}{2!}\frac{x}{2}+\dots+\frac{1}{r!}(\frac{x}{2})^r\dots)f ]+i[b(\alpha\beta-\gamma\delta)e^q (1-\frac{1}{2!}\frac{x}{2}+\dots+\\&+(-1)^r\frac{1}{r!}(\frac{x}{2})^r\dots)f]+ic\gamma f+i^2d\partial_xf+\frac{1}{2}d\partial_qf+ib\gamma e^q[i\partial_xf -\frac{1}{2i}\mathcal F_p\left(pi\mathcal F^{-1}_p(xf)\right)+\\& +\dots+(\frac{1}{2i})^r\frac{1}{r!}(-1)^r\mathcal F_p\left(pi^r\mathcal F^{-1}_p(x^rf)\right)+\dots]=d(\frac{1}{2}\partial_q-\partial_x)f+\\&+[iae^{-(q-\frac{x}{2})}+ib(\alpha\beta-\gamma\delta)e^{q-\frac{x}{2}}]f+ic\gamma f+\frac{1}{2}e^{-\frac{x}{2}}b\gamma e^q\partial_qf-b\gamma e^qe^{-\frac{x}{2}}\partial_xf\\&=(d+ b\gamma e^{q-\frac{x}{2}})(\frac{1}{2}\partial_q-\partial_x)f+[iae^{-(q-\frac{x}{2})}+ib(\alpha\beta-\gamma\delta)e^{q-\frac{x}{2}}+ic\gamma]f
\end{array}$$
\item At last, if $\tilde A$ is defined by (10) then :
$$\hat\ell_A(f)=(d- b\gamma e^{q-\frac{x}{2}})(\frac{1}{2}\partial_q-\partial_x)f+[-iae^{-(q-\frac{x}{2})}-ib(\alpha\beta-\gamma\delta)e^{q-\frac{x}{2}}+ic\gamma]f$$
\end{enumerate}
The theorem is therefore proved.\hfill$\square$

\par
\begin{remark}{\rm
Setting new variables $s = q - \frac{x}{2}$, $t = q + \frac{x}{2}$, we have 
$$\hat{\ell} _A(f)=\cases$$\big(d\partial_s+ia\alpha e^{-s}\big)f\vert_{(s,t)} $$ &\text{if $\tilde {A}$  is defined by (6) }\\
$$\big(d\partial_s+ib\beta e^s\big)f\vert_{(s,t)} $$ &\text{if $\tilde {A}$ is defined by  (7)}\\
$$\big(d\partial_s +i[a\alpha e^{-s}+ b\beta e^s]\big)f\vert_{(s,t)} $$ &\text{if $\tilde {A}$ is defined by (8)}$$ \\
$$\Big((d+b\gamma e^s)\partial_s+$$\\$$i[ae^{-s}+b(\alpha\beta-\gamma\delta)e^s +c\gamma ]\Big)f\vert_{(s,t)} .&\text{if $\tilde {A}$ is defined by  (9)$$}\\
$$\Big((d-b\gamma e^s)\partial_s+$$\\$$i[-ae^{-s}-b(\alpha\beta-\gamma\delta)e^s +c\gamma ]\Big)f\vert_{(s,t)} .&\text{if $\tilde {A}$ is defined by  (10)$$} 
\endcases$$
}\end{remark}

\begin{theorem}
With above notations we obtain the operators : $$\hat{\ell} _A=\cases$$\hat\ell_A^{(2)}=\big(d\partial_s+ia\alpha e^{-s}\big)\vert_{(s,t)} $$ &\\
$$\hat\ell_A^{(3)}=\big(d\partial_s+ib\beta e^s\big)\vert_{(s,t)} $$ &\\
$$\hat\ell_A^{(4)}=\big(d\partial_s+i[a\alpha e^{-s}+ b\beta e^s]\big)\vert_{(s,t)} $$ & \\
$$\hat\ell_A^{(5)}=\Big((d+b\gamma e^s)\partial_s+i[ae^{-s}+b(\alpha\beta-\gamma\delta)e^s+c\gamma]\Big)\vert_{(s,t)} .&  \\
$$\hat\ell_A^{(5')}=\Big((d-b\gamma e^s)\partial_s+i[-ae^{-s}-b(\alpha\beta-\gamma\delta)e^s+c\gamma]\Big)\vert_{(s,t)} 
\endcases$$ which provides the representations of the Lie algebra $\mathfrak g$=$Lie(\mathbb R\ltimes \mathbb H_3)$.\par
Furthermore, $\forall A,B \in\mathfrak g$,
$$\hat\ell_A\circ\hat\ell_ B- \hat\ell_B\circ\hat\ell_ A =\hat\ell_{[A, B]}$$
\end{theorem}

{\it Proof}    
For each compactly supported $C^\infty$ function $f \in C^\infty_0(\mathbb R^2)$ and for  $A,B \in Lie(\mathbb R \ltimes \mathbb H_3)$,we have 
$$\hat\ell_{(\mu_1A+\mu_2 B)}(f)=\mathcal F_p \circ \ell_{(\mu_1A+\mu_2 B)}\circ \mathcal F^{-1}_p(f)=\mathcal F_p\Big(i(\widetilde{\mu_1A+\mu_2B})\star \mathcal F^{-1}_p\Big)=$$
$$=\mu_1\mathcal F_p\circ \ell_A\circ\mathcal F^{-1}_p(f)+\mu_2\mathcal F_p\circ \ell_B\circ\mathcal F^{-1}_p(f)=\mu_1\hat\ell_A(f)+\mu_2\hat\ell_B(f)\quad\forall \mu_1,\mu_2\in\mathbb R. $$
Moreover,
$$\hat\ell_A\circ\hat\ell_B(f)-\hat\ell_B\circ\hat\ell_A(f)=\hat\ell_A \Big(\mathcal F_p \circ \ell_B \circ \mathcal F^{-1}_p(f)\Big)- \hat\ell_B \Big(\mathcal F_p \circ \ell_A \circ \mathcal F^{-1}_p(f)\Big)=$$
$$=\mathcal F_p\Big(i\tilde A\star(i\tilde B\star\mathcal F^{-1}_p(f)\Big)-\mathcal F_p \Big(i\tilde B\star (i\tilde A\star\mathcal F^{-1}_p(f)\Big)=\mathcal F_p\Big(i\widetilde{[A,B]}\star\mathcal F^{-1}_p(f)\Big)=\hat\ell_{[A,B]}(f)$$\hfill$\square$

\par
\begin{definition}{\rm
Let $\Omega_F^{\lambda}$  be K-orbits of the real diamond Lie group  $G$. With  $A$ runs over the Lie algebra $\mathfrak g=Lie(G)$, 
\begin{itemize}
\item $(\Omega^2 ,\hat\ell_A^{(2)}) ;(\Omega^3 ,\hat\ell_A^{(3)})$ are called the  quantum half-planes,
\item $(\Omega^4 ,\hat\ell_A^{(4)})$  - quantum hyperbolic   cylinder,
\item  $(\Omega^5 ,\hat\ell_A^{(5)} ,\hat\ell_A^{(5')})$ - quantum hyperbolic paraboloid, 
\end{itemize}
with respect to the co-adjoint action of  Lie group $G$. In the other words, $(\Omega_F,\hat\ell_A)$, with $A$ running over the Lie algebra $\mathfrak g$ is called a {\it quantum co-adjoint orbit} of  Lie group $G$.  
}\end{definition}
As G=$\mathbb R\ltimes \mathbb H_3$   is connected and simply connected, we obtain a unitary representations \quad $ T$ of G defined by the following formula $$ T(\exp A) :=\exp(\hat{\ell}_A);\quad  A\in\mathfrak g$$  More detail, $$\exp(\hat{\ell}_A) =\cases$$\exp (d\partial_s+ia\alpha e^{-s})\vert_{(s,t)} $$ &\text{if $\tilde {A}$ is defined by  (6) }\\
$$\exp (d\partial_s+ib\beta e^s)\vert_{(s,t)} $$ &\text{if $\tilde {A}$  is defined by (7)}\\
$$\exp (d\partial_s+i[a\alpha e^{-s}+ b\beta e^s])\vert_{(s,t)} $$ &\text{if $\tilde {A}$  is defined by (8)}$$ \\
$$\exp ((d+b\gamma e^s)\partial_s+$$\\$$i[ae^{-s}+b(\alpha\beta-\gamma\delta)e^s+c\gamma]   )\vert_{(s,t)}.&\text{if $\tilde {A}$  is defined by (9)$$}\\ 
$$\exp ((d-b\gamma e^s)\partial_s+$$\\$$i[-ae^{-s}-b(\alpha\beta-\gamma\delta)e^s+c\gamma]   )\vert_{(s,t)}.&\text{if $\tilde {A}$  is defined by (10)$$}
\endcases$$

This means that we refind all the representations $ T(\exp A)$ of the real diamond Lie group $\mathbb R\ltimes \mathbb H_3$, those could implicitly obtained from (induction) orbit method induction. What we did here gives us more precise analytic formulas in this case for orbit method induction.

\vfil\eject

$$\text{ ACKNOWLEDGMENT}$$
The author would like to express his gratitude to Professor Do Ngoc Diep  for all his helpfulness and for suggesting many of the topics considered in this paper. The author also thanks Dr. Nguyen Viet Dung for his encouragement.

\end{document}